\documentclass[12pt]{article}
\usepackage{amsmath,amsfonts,amssymb,fullpage}

\def\bull{\vrule height .9ex width .8ex depth -.1ex}
\begin{document}

\title{Notes on the Second Eigenvalue of the Google Matrix}
\author{Roger Nussbaum\thanks{Partially supported by NSF
DMS-00-70829\newline
}}

\date{\today}
\maketitle

\begin{abstract}
If $A$ is an $n\times n$ matrix whose $n$ eigenvalues are ordered in
terms of decreasing modules, $|\lambda_1 | \geq |\lambda_2| \geq
\cdots |\lambda_n|$, it is often of interest to estimate
$\frac{|\lambda_2|}{|\lambda_1|}$.  If $A$ is a row stochastic matrix
(so $\lambda_1 = 1$), one can use an old formula of R. L. Dobrushin to
give a useful, explicit formula for $|\lambda_2|$.  The purpose of
this note is to disseminate these known results more widely and to
show how they imply, as a very special case, some recent theorems of
Haveliwala and Kamvar about the second eigenvalue of the Google matrix.
\end{abstract}



If $A=(a_{ij})$ is an $n\times n$ real matrix, $A$ has $n$ (counting
algebraic multiplicity) complex eigenvalues which can be listed in
order of decreasing  modules: $|\lambda_1|\geq |\lambda_2|\geq \cdots
\geq |\lambda_n|$.  We have $|\lambda_1|=\sup \{
|\lambda|\big|\lambda$ is an eigenvalue of $A\}$ and $|\lambda_1|$ is
called the spectral radius of $A, r(A): = |\lambda_1|$.  In many
problems it is of interest to estimate
$\frac{|\lambda_2|}{|\lambda_1|} = \frac{|\lambda_2|}{r(A)}$.  Indeed,
an analogous problem is of great interest for bounded linear maps on
Banach spaces: {see} [2], [3] and the references there.

Slightly more generally, suppose that $V$ is an $m$-dimensional real
vector space and $L:V\to V$ is a linear map.  Again $L$ has $m$
possibly complex eigenvalues which can be written in order of
decreasing modules: $|\lambda_1| \geq |\lambda_2| \geq \cdots \geq
|\lambda_m|$.  If $\|\cdot\|$ denotes any norm on $V$ (recall that all
norms on a finite dimensional real vector space give the same
topology), we can define

\begin{equation}
\|L\| = \sup\{\|Ly\|: y\in V, \|y\|\leq 1\}.
\end{equation}
It is known that 
\begin{eqnarray}
r(L) = |\lambda_1| &=& \sup\{|\lambda|:\lambda\text{ is an eigenvalue of }L\}\\
                &=&\underset{k\to \infty}{\lim }\| L^k\|^{\frac{1}{k}}=
\underset{k\geq 1}{\inf} \| L^k\|^{\frac{1}{k}},\nonumber
\end{eqnarray}
where $L^k$ denotes the composition of $L$ with itself $k$-times.

We shall consider elements of $\Bbb R^n$, as usual, as column vectors.
An $n\times n$ matrix $B$ induces a linear map $\Lambda : \Bbb R^n \to
\Bbb R^n$ by $\Lambda(y) = By$.  If $V$ is a vector subspace of $\Bbb
R^n$ and $By\in V$ for all $y\in V$, then $B$ induces a linear map
$L:V\to V$ by $L(y) = By$ for $y \in V$.  If $\dim (V) = m$, then $L$
has (counting algebraic multiplicity) precisely $m$ eigenvalues, and
these are the eigenvalues of $B$ whose corresponding eigenvectors lie
in the complexification of $V$.

Now suppose that $A = (a_{ij})$ is an $n\times n$ row stochastic
matrix, so $a_{ij} \geq 0$ for all $i, j$ and $\sum^n_{j=1} a_{ij} = 1 $ for
$ 1\leq i \leq n$.  Denote by $x^t$ the transpose of a vector $x$ and
by $B^t$ the transpose of a matrix $B$.  If $e=(1, 1, \cdots, 1)^t$,
then $Ae=e$, so $1\in \sigma (A)$, where $\sigma(A)$, the spectrum of
$A$, denotes the collection of eigenvalues of $A$.  Recall that (in
general) $\sigma (A) = \sigma (A^t)$, so $1\in \sigma (A^t)$.  It follows
that (in general) $r(A) = r (A^t)$; and it is an elementary fact (the
proof is sketched below) that $r (A) = 1$ for $A$ row stochastic.

It will be convenient to use the $L^1$ norm $\| \cdot\|$, on $\Bbb
R^n$, so for 
\begin{eqnarray}
y &=& (y_1, y_2, \cdots, y_n)^t \in \Bbb R^n\\
              \|y \|_1: & =& \sum^n_{i=1} |y_i|  \nonumber.
\end{eqnarray}
Using the $L^1$ norm, we get a corresponding norm on $n\times n$
matrices $B=(b_{ij})$, since these matrices induce linear maps:

\begin{equation}
\|B\|_1 = \sup\{ \| By\|_1 : \|y\|_1 \leq 1, y\in \Bbb R^n\}.
\end{equation}
Indeed, using this norm when $A$ is row stochastic, it is easy to see
that $\|A^t\|_1=1$.  Since $1\in \sigma (A^t)$, we deduce, using
eq. (2), that $r (A^t)=1$ and hence $r(A) = 1$.

\begin{equation}
\text{ If } x,y \in \Bbb R^n, \text{ let } \langle x,y\rangle =
\sum^n_{i=1} x_i y_i. \nonumber
\end{equation}

\newtheorem{lem}{Lemma}\begin{lem}  Let $A$ be an $n\times n$ row
stochastic matrix and let $V=\{x = (x_1, x_2, \cdots, x_n)^t\in \Bbb
R^n| \sum^n_{1} x_i = 0\}$.  If $x\in V$, it follows that $A^t x\in V$.
\end{lem}

\noindent{\it Proof:}  If $x\in V, \langle x,e\rangle
= 0 =\langle x, Ae\rangle = \langle A^t x, e\rangle$, so $A^t x \in
V  $.  $\quad \bull$
\medskip

Henceforth, $V$ will be as in Lemma 1.

If $A$ is row stochastic, let $L: V\to V$ be the linear map induced by
$A^t$.  If $1=\lambda_1 \geq |\lambda_2|\geq \cdots \geq |\lambda_n|$
are the moduli of the eigenvalues of $A^t$, our previous remarks show
that $\lambda_2, \lambda_3, \cdots, \lambda_n$ are the eigenvalues of
$L$ and eq. (2) implies that

\begin{equation}
|\lambda_2|:= \text{ spectral radius of } L: = r(L) = \underset{k\to
\infty}{\lim} \|L^k\|_1^{\frac{1}{k}} = \underset{k\geq 1}{\inf}\|L^k\|_1^{\frac{1}{k}}.
\end{equation}

By definition,
\begin{eqnarray}
\|L^k\|_1 &=& \sup\{\|(A^k)^t y \|_1: y \in V, \|y \|_1\leq 1\}\\
          &=& q(A^k). \nonumber
\end{eqnarray}
Note that $A^k$ is a row stochastic matrix.  If $B$ is any row
stochastic matrix, we follow (6) and define

\begin{eqnarray}
q(B) &=&\sup\{\|B^t y\|_1: y\in V, \|y\|_1\leq 1\},\text{ where }\\
   V &=& \{y\in \Bbb R^n : \sum^n_{1} y_i = 0\} \nonumber.
\end{eqnarray}
The formula given by eqns (5)-(7) would be of limited usefulness
without an explicit formula for $q(B)$.  Fortunately, Dobrushin has
given such a formula in Lemma 1, Section 3 of [1]; a slightly more
general result is proved in Lemma 3.4 of [5].

\begin{lem} (Dobrushin [1]). Let
$B=(b_{ij})$ be an $n\times n$ row stochastic matrix and let $V=\{y\in
\Bbb R^n|\sum^n_1 y_i = 0\}$.  If $q(B)$ is defined by (7), then
\end{lem}

\begin{equation}
q(B) = \left(\frac{1}{2}\right)\underset{i,k}{\sup} \left(\sum^n_{j=1}
|b_{ij} - b_{kj}|\right) = 1-\underset{i,k}{\min}\sum^n_{j=1} \min (b_{ij}, b_{kj}).
\end{equation}

Combining the above observations we obtain a useful formula for
$|\lambda_2|$ when $A$ is row stochastic.

\newtheorem{theo}{Theorem}\begin{theo}  Let $A$ be an $n\times n$ row
stochastic matrix with eigenvalues $1=\lambda_1, \lambda_2, \cdots,
\lambda_n$, where  $1= |\lambda_1|\geq |\lambda_2| \geq \cdots
\geq |\lambda_n|$ and eigenvalues are counted with algebraic
multiplicity. (Recall that these eigenvalues are the same as the
eigenvalues of $A^t$).  Then we have 
\end{theo}
\begin{equation}
|\lambda_2|= \underset{k\to\infty}{\lim} q(A^k)^{\frac{1}{k}} =
\underset{k\geq 1}{\inf} q(A^k)^{\frac{1}{k}},
\end{equation}
where $q(B)$ is defined by eq. (8)

One can verify directly that for {\bf any} $n\times n$ real matrix
$B$, if we {\bf define} $q(B): = (\frac{1}{2})\sup$ $\left(\sum^n_{j=1}
|b_{ij} - b_{kj}|\right)$, then $q$ is a seminorm, ie, $q(B+C) \leq
q(B) + q(C)$ for any $n\times n$ matrices $B$ and $C$ and $q(\alpha B)
= |\alpha| q (B)$ for any scalar $\alpha$.  Furthermore, if $B$ and
$C$ are any $n\times n$ real matrices which have $e^t = (1, 1, \cdots,
1)^t$ as an eigenvector, then $q(BC) \leq q (B) q (C)$.

For the case that $E$ is row stochastic  of rank 1, the next result is
proved in [4] by different methods.

\newtheorem{cor}{Corollary} \begin{cor}Let $P$ be an $n\times n$ row stochastic
matrix and let $E$ be an $n\times n$ row stochastic matrix.  If $0\leq
c\leq 1,$ let $A= cP + (1-c) E$ and let $\lambda_2$ denote the second
eigenvalue of $A$ (where $1=\lambda_1 \geq  |\lambda_2 | \geq
|\lambda_3|\geq \cdots \geq |\lambda_n|$).  Then we have 
\end{cor}

\begin{equation}
|\lambda_2| \leq c q (P) + (1-c) q(E)\text{ and } |\lambda_2|\leq
c\text{ if } E \text{ has rank } 1.
\end{equation}

\noindent{\it Proof:} Theorem 1 implies that
$|\lambda_2| \leq q(A) \leq cq(P) + (1-c)q(E)$.  If $E$ has rank 1,
then all rows of $E$ are identical and $q(E)=0$.  Since $q(P) \leq
1$, this gives $q(A)\leq c$  when $E$ has rank 1. $ \quad \bull$

The following result is proved in [4] by a somewhat more involved
argument.

\begin{cor}
Let $P$ and $E$ be as in Corollary 1 and assume $E$ has rank 1.  If
$\{y\in \Bbb R^n|y = Py\}$ has dimension greater than one then
$|\lambda_2|=c$, where $\lambda_2$ is the second eigenvalue of $c P
+(1-c)E$.  In fact, we have $\lambda_2=c$.
\end{cor}

\noindent{\it Proof:}  By assumption, there are
linearly independent vectors $v$ and $w$ with $Pv=v$ and $Pw=w$.  We
also know there is a vector $z\in \Bbb R^n, \sum^n_{i=1} z_i = 1,
z_i\geq 0$ for $1\leq i \leq n$, such that $Ex = \langle x, z\rangle
e$.  If $\langle v, z\rangle = 0, (c P + (1-c) E) v = c Pv = cv,$ and
if $\langle w, z\rangle = 0, (cP+(1-c) E)w=cw$.  Thus assume that
$\langle v, z\rangle\neq 0$ and $\langle w, z\rangle \neq 0$ and
define $\xi =-\langle w, z\rangle v+ \langle v, z\rangle w$.  Note
that $\langle \xi, z\rangle = 0$ and $\xi \neq 0$ because $v$ and $w$
are linearly independent.  It follows that $(cP + (1-c)E)\xi = cP\xi = c\xi$.  We conclude that $c$ is an eigenvalue of $cP+(1-c) E$.
Since we already know from Corollary 1 that $|\lambda_2|\leq c$ we
conclude that $|\lambda_2|=c$ and that we can take $\lambda_2=c$. $\bull$
\medskip

\noindent{\bf Remark 1 }  The statement (see [4]) that
``$P$ has at least two irreducible closed subsets'' is equivalent to
the assertion that $\dim\{ y\in \Bbb R^n |y = Py\} \geq 2$.
\medskip

\noindent{\bf Remark 2 }  Suppose that $A$ is row
stochastic and $q(A) = \kappa <1$.  Let $\sum = \{x \in \Bbb R^n|\sum^n_1
x_i = 1$ and $x_j \geq 0$ for $1\leq j \leq n\}$.  One easily checks
that $A^t x \in \sum$ if $x\in \sum$, and our previous remarks show
that $\|A^t x - A^t y \|_1 \leq \kappa \| x-y \|_1$ for all $x, y \in
\sum$.  By the contraction mapping theorem, for any $x\in \sum,
\underset{k\to \infty}{\lim} (A^t)^k x=v$, where $A^t v=v$.  Also, the
rate of convergence can be estimated in terms of $\kappa$.  The same
assertions are true if $q(A^m)=\kappa_m < 1$ for some $m\geq 1$.

\noindent Mathematics Department, Hill Center\\
Rutgers University\\
110 Frelinghuysen Road\\
Piscataway, New Jersey\\
U.S.A. 08854-8019

\end{document}